\documentstyle[12pt]{article}
\begin{document}
\title{
Poincare'-Birkhoff-Witt property for bicovariant differential algebras
on simple quantum groups.
}
\author{G.E.Arutuynov\thanks{Steklov Mathematical Institute, Vavilov
42, GSP-1, 117966,
Moscow,
Russia;$~~~~~~~~~~~~$ $~~~~~~$arut@qft.mian.su},$~~$
A.P.Isaev\thanks{Laboratory of Theoretical Physics, JINR,
Dubna SU-101 000 Moscow,
Russia; $~~~~~~~~~~$ $~~~~~$
isaevap@thsun1.jinr.dubna.su}
\\
and\\
Z.Popowicz\thanks
{Institute of Theoretical Physics, University of Wroclaw, Pl.M.Borna
9,
50-204 Wroclaw, Poland; ziemek@ift.uni.wroc.pl}
}
\date {December 1994~}
\maketitle
\begin{abstract}
We investigate the possibility to construct
bicovariant differential calculi on quantum groups $SO_q(N)$ and
$Sp_q(N)$ as a quantization of an underlying bicovariant bracket.
We show that, opposite to $GL(N)$ and
$SL(N)$-cases, neither of possible graded
$SO$- and $Sp$- bicovariant brackets (associated
with a quasitriangular $r$-matrices)
obey the Jacobi identity when the differential forms are
Lie algebra-valued.
The absence of a classical Poisson structure gives an
indication that differential algebras describing bicovariant
differential
calculi on quantum orthogonal and symplectic groups are not of
Poincar\'{e}-Birkhoff-Witt type.
\end{abstract}
\newpage
Bicovariant differential calculus (BDC) on quantum groups
initiated by Worono\-wicz's \cite{Wor} provides a meaningful example
of the
noncommutative differential geometry \cite{Conn}.
On the other hand, it serves as the starting point to formulate a new
class of gauge theories with a simple quantum group playing the role
of a
gauge group \cite{AV},\cite{IsPop}. Many of the phenomena, which
one can encounter studying these theories, are not privileges of a
model
but have their origin in the theory of BDC. Thus, it is extremely
important
to investigate the general properties of BDC on simple quantum groups.

In this letter we are aimed to find out whether external algebras
on quantum groups $SO_q(N)$ and $Sp_q(N)$ are of Poincare'-Birkhoff-
Witt
(PBW) type, {\em i.e.} they possess a unique basis of
lexicographycally ordered monomials.
This is not an academic question, since it has a strong influence
on all differential geometry associated with these groups. In
particular,
the missing of PBW property under quantization means that the
classical
and corresponding to it quantum system
has different number of observables.

Our consideration is based on the $R$-matrix approach of \cite{Fad},
which is very useful in dealing with BDC's.

Recall that the central point of Woronowicz's theory is the
constructing of bicovariant bimodules $\Gamma$ over a Hopf algebra
${\cal
A}$ (the algebra of functions
on a quantum group). The bimodules over
$\cal A$ supplied with two coactions:
\begin{equation}
\Delta_R:~\Gamma\rightarrow \Gamma\otimes {\cal A}~~~\mbox{and}~~~
\Delta_L:~\Gamma\rightarrow {\cal A}\otimes \Gamma
\label{coact}
\end{equation}
satisfying the set of axioms \cite{Wor}. Bicovariant bimodules
are interpreted as noncommutative analogues of tensor bundles
$\Gamma_{cl}$
over Lie groups. For the case of simple quantum groups the
classification of bicovariant bimodules was obtained in \cite{J} and
confirmed in \cite{Schmudgen}.

A first order differential calculus is defined as a pair
$(\Gamma, {\bf d})$, where
differential ${\bf d}:~{\cal A}\rightarrow\Gamma$ is a
nilpotent mapping obeying the Leibnitz rule.

The bicovariant wedge product of
two left-invariant 1-forms $\Omega_{i}$ is defined via tensor
algebra construction:
\begin{equation}
\label{w6}
\begin{array}{c}
\Omega_{i} \wedge \Omega_{j} =
\Omega_{i} \, \otimes_{\cal A} \, \Omega_{j} - \sigma^{lk}_{ij}
\Omega_{l} \, \otimes_{\cal A} \, \Omega_{k} \,
\end{array}
\end{equation}
or in the concise matrix notation:
\begin{equation}
\label{w7}
\Omega_{1} \wedge \Omega_{2} =
( I_{12} - \sigma_{12} ) \Omega_{1} \, \otimes_{\cal A} \, \Omega_{2}
\; .
\end{equation}
Here matrix $\sigma_{12}$ satisfies the Yang-Baxter equation (YBE)
\cite{Wor}:
\begin{equation}
\label{w8}
\sigma_{23} \sigma_{12} \sigma_{23} =
\sigma_{12} \sigma_{23} \sigma_{12} \; ,
\end{equation}
which provides the assosiativity of the wedge product:
\begin{equation}
\label{w7a}
(\Omega_{1} \wedge \Omega_{2}) \wedge \Omega_{3}  =
\Omega_{1} \wedge ( \Omega_{2} \wedge \Omega_{3} )  \; .
\end{equation}

Therefore, adopting (\ref{w6}),(\ref{w7}) one can construct,
starting from
$\Gamma$,  an associative external algebra
$\Gamma^{\wedge}=\sum_{n} \Gamma^{(n)}$,
where $\Gamma^{(0)} = {\cal A}, \; \Gamma^{(1)} = \Gamma$ and
$\Gamma^{(n)}$ is the space of $n$-forms. It is proved \cite{Wor} that
a first order differential calculus can be lifted to higher order
differential forms via extending $\Gamma$ by an additional bi-
invariant
1-form $X$ generating ${\bf d}:~~{\bf d}\Omega=[X,\Omega]_{\pm}$,
$\Omega\in \Gamma^{\wedge}$.

However, let us stress that YBE (\ref{w8}) does not
guarantee that $\Gamma^{\wedge}$ constructed in such a way is the
algebra of PBW type. From the point of view of
general BDC theory \cite{Wor} the fulfillment
of PBW property for quantum external algebras is an additional
physical requirement.

The general properties of $\sigma$ for quantum simple Lie groups
(in particular, the projector expansion)
were studied in \cite{Wat},
\cite{Cast}. In \cite{Wat} the authors modified the
definition (\ref{w6}) by imposing the additional quadratic relations
on the generators $\Omega_i$ (we will comment this later).
The direct investigation of the PBW property for (\ref{w6}) is rather
involved and to our knowledge it has been tackled only for the
external algebras on $GL_q(N)$ and $SL_q(N)$ \cite{Man},\cite{IP},
\cite{FP}.
Fortunately, for quantum groups
$SO_q(N)$ and $Sp_q(N)$ $(N=2n)$, one can
exploit ideas of \cite{SKL}, \cite{SEM} that the quantum
groups can be obtained by quantizing of classical Poisson structures
and try to answer the question, about PBW property, just on the
semiclassical level. This is due to the existence of the
infinitesimal version of bicovariant differential calculi on the
quantum
groups provided by graded bicovariant brackets \cite{AM}
which can be deduced directly as a semiclassical
limit of the PBW algebras presented in \cite{IP},\cite{FP},
\cite{Zum}. In this approach,
it is assumed that
the anticommutator in $\Gamma$ is determined,
in the semiclassical approximation,
by a graded bracket on
$\Gamma_{cl}$:
\begin{equation}
\Omega\wedge\Omega'+\Omega'\wedge\Omega=\hbar\{\Omega,\Omega'\}+
\hbar^2(\ldots),
\label{cl}
\end{equation}
satisfying the condition of bicovariance:
\begin{equation}
\begin{array}{l}
\Delta_{L,R}(\{ \Omega, \; \Omega'\})
=\{\Delta_{L,R}( \Omega ) , \, \Delta_{L,R} (\Omega') \},~~~
\Omega,\Omega' \in \Gamma_{cl}^{1}.
\end{array}
\label{y4}
\end{equation}
The coactions $\Delta_L: \Gamma \rightarrow {\cal A}\otimes\Gamma$
and $\Delta_R: \Gamma \rightarrow \Gamma \otimes {\cal A}$
coming in (\ref{y4}) are the classical analogues of (\ref{coact}).
On the matrix elements $\Omega_{i}^{j}$ ($i,j = 1,2, \dots N$)
of the Lie-valued
left-invariant Cartan's
form $\Omega$ ($\Omega^{i}_{j}$ generates $\Gamma_{cl}^{\wedge}$)
$\Delta_{L,R}$ are defined as
\begin{equation}
\Delta_{L}(\Omega^{i}_{j}) = I \otimes \Omega^{i}_{j} \equiv
\Omega^{i}_{j},
{}~~~ \Delta_{R}(\Omega^{i}_{j}) = \Omega^{k}_{l} \otimes
(T^{-1})^{i}_{k} T^{l}_{j} \equiv (T^{-1} \, \Omega \, T)^{i}_{j},
\label{y3}
\end{equation}
and to arbitrary element of $\Gamma_{cl}^{\wedge}$ they are
extended as homomorphisms with respect to the wedge product.
We introduce the matrix element $T^{i}_{j}$ of the groups
in (\ref{y3}) and, for the sake of simplicity,
we omit the signs of the tensor
products in the last parts of (\ref{y3}).
In the following we will refer to bracket satisfying (\ref{y4})
as to bicovariant with respect to the $L,R$-coactions of $G$.

To turn $\Gamma_{cl}^{\wedge}$ into a classical Poisson system we also
impose on the graded bracket in (\ref{cl}) the following conventional
requirements :

i) the symmetry condition:
$$
\{\rho, \, \rho' \} = (-1)^{\deg{(\rho)}\deg{(\rho')}+1}
\{ \rho', \, \rho \} \; ,
$$

ii) the graded Jacobi identity:
\begin{equation}
\label{y5}
(-1)^{\deg{\rho_1}\deg{\rho_3}}
\{\{ \rho_1 , \, \rho_{2} \}, \, \rho_3 \} +
(-1)^{\deg{\rho_2}  \deg{\rho_3} }
\{\{ \rho_3 , \, \rho_{1} \}, \, \rho_2 \} +
\end{equation}
$$
(-1)^{\deg{\rho_1}  \deg{\rho_3}}
\{\{ \rho_2 , \, \rho_{3} \},
\, \rho_1 \} = 0 \; ,
$$
and

iii)
$$
\{ \rho_{1} \otimes \rho_{2}, \, \rho_{3} \otimes \rho_{4} \} =
\{ \rho_{1} , \, \rho_{3} \} \otimes \rho_{2} \, \rho_{4}  +
\rho_{1} \, \rho_{3} \otimes \{ \rho_{2}, \, \rho_{4} \} \; .
$$
In addition we demand (as usual) this bracket to be a graded
differentiation:
$$
\{ \rho_1 , \, \rho_2 \wedge \rho_3 \} =
\{ \rho_1 , \, \rho_2 \} \wedge \rho_3 +
(-1)^{\deg{\rho_1} \cdot \deg{\rho_2}}
\rho_{2} \wedge \{ \rho_1 , \, \rho_{3} \}
$$
If a bracket satisfying the requirements above exists,
then $\Gamma_{cl}^{\wedge}$ is said to be equipped with a graded
Poisson-Lie (PL) structure \cite{AM} and one can consider
$\Gamma^{\wedge}_{cl}$ as a phase space for some
graded dynamical system.

Thus, in general, an algebra of quantum external
forms is expected to be a graded bicovariant algebra with the graded
commutator that produces in the semiclassical limit a graded
bicovariant bracket. Now the fact that this algebra is of PBW type
leads, in semiclassics, to the
requirement on the corresponding bicovariant
bracket to be Poisson, {\em i.e.}, to satisfy the Jacobi identity
(here and below we confine ourselves only
with a consideration of the exterior
algebras (\ref{w6}) having usual classical limit).

In the cases of $GL(N)$ and $SL(N)$, graded PL structures exist
\cite{AM}, \cite{AAM} and the corresponding
algebras of quantum external forms
are of PBW type \cite{IP},\cite{FP}. Moreover, if a graded PL
structure exists,
then bicovariance and PBW property can be considered as main
quantization principles.

Let $G$ be $SO(N)$ or $Sp(N)$ groups and $\cal G$ be the corresponding
Lie algebra. The following terminology (see \cite{DR}) will be
useful. A skewsymmetric solution $r$ $(r\in {\cal G}\wedge {\cal G})$
of the classical YBE (CYBE) will be refered as a triangular $r$-matrix
and a skewsymmetric $r$ obeying the modified YBE (mYBE) will be
refered
as a quasitriangular one.

Our strategy is as follows. It is natural to consider a graded
PL structure on $\Gamma^{\wedge}_{cl}$ generated by the brackets
between the components of $\cal G$-valued Cartan's form $\Omega$,
i.e., when $\Omega^{i}_{j} = \omega_{\alpha} (t^{\alpha})^{i}_{j}$
where $t^{\alpha} \in {\cal G}$.
It is shown below
(see also \cite{AR}) that these brackets are defined via a triangular
$r$-matrix. However, this is not for the case of the standard $r$-
matrix
\cite{DR} assosiated with simple Lie algebras. In other words,
if we employ a quasitriangular $r$ (that is relevant for subsequent
quantization) and  require $G$-covariant Poisson bracket
$\{\Omega_1,\Omega_2\}$ to be an element of
${\cal G}\otimes {\cal G}$ (as a matrix), then we get a unique
solution $\{\Omega_1,\Omega_2\}=0$ (see below).
One may hope that by discarding the requirement $\Omega\in {\cal G}$
it would be possible to obtain a graded bicovariant bracket with a
quasitriangular $r$. Below we analyze this possibility and, hence,
assume the general situation when $\Omega\in Mat(N,{\bf C}) \sim $
$ gl(N,{\bf C})$.
In this case, a covariancy group $G$ of graded brackets on $gl(N)$
is a subgroup of $GL(N)$.
Let us stress that this is in agreement with the quantum external
algebra
construction \cite{Wor} where $\dim{\Gamma^{(1)}}=N^2$.
One remark is in order.
We will not consider in this letter
the graded bicovariant brackets which are covariant under the groups
isomorphic with the
linear groups of $A_{n-1}$ series
(e.g. $SO(3) \sim$ $Sp(2) \sim$ $SL(2)$, see \cite{AR} for
discussion).

Now we recall briefly the basic facts about Lie groups $G$
corresponding
to $so(N)$ or $sp(N) ~(N=2n)$ Lie algebras. The fundamental
representation
of $G$ is given by
$$
TCT^tC^{-1}=CT^{t}C^{-1}T=I,
$$
where $N\times N$ metric $C$ is $C^{ij}=\delta^{ij'}$ for $SO(N)$ and
$C^{ij}=\epsilon_i\delta^{ij'}$ for $Sp(N)$, $i'=N+1-i$,
$\epsilon_i=1$
$(i=1,\ldots, n)$, and $\epsilon_i=-1$ if $(i=n,\ldots, 2n)$.
We denote by $C^{ij}$ ($C_{ij}$) matrix elements of $C$ $(C^{-1})$.

The fundamental representations of the corresponding Lie algebras
are defined as follows
$$
{\cal G}=\{X\in Mat(N,{\bf C})|~X^t=-CXC^{-1}\}.
$$

To simplify the calculations we introduce an operation $\tilde{}$
acting
on $\Omega,\Omega^2$, etc., in the following way
\begin{equation}
\tilde{\Omega}=C\Omega^t C^{-1},~~~
\widetilde{\Omega^2}=C({\Omega^2})^t C^{-1}=
-({\tilde{\Omega}})^2.
\label{til}
\end{equation}
Clearly, $\tilde{\tilde{\Omega}}=\Omega$.
Using this operation we split matrix-valued forms as
$\Omega^{\pm}=\Omega \pm \tilde{\Omega}$.
Note that the form $\Omega^{-}$ belongs to $\cal G$ in the
fundamental representation.

It can be shown \cite{AM1} that the general form of a $Z$-graded
bicovariant bracket $\{ \Omega_{1}, \, \Omega_{2} \}$ is
\begin{equation}
\label{y12}
\{ \Omega_{1}, \; \Omega_{2} \} =
\left[ \Omega_{1}, \; [ \Omega_{2} , \; r_{12} ] \right]_{+} +
Tr_{34}(W_{1234} \Omega_{3} \Omega_{4} ) \; ,
\end{equation}
where $r_{12}$ is the quasitriangular $r$-matrix and
$W_{1234}$
is a $G$-invariant tensor:
\begin{equation}
\label{y10}
W_{1234} = T_1 \, T_2 \, T_3 \, T_4 \, W_{1234}
T^{-1}_{1} T^{-1}_{2} T^{-1}_{3} T^{-1}_{4} \; ,
\end{equation}
with symmetry properties:
\begin{equation}
\label{y11}
W_{1234} = W_{2134} = -W_{1243} \; ,
\end{equation}
Here indeces $1,2,3,4$ denote the numbers
of the matrix spaces.
Thus, to construct a
general $SO$- $(Sp)$-bicovariant bracket we have to enumerate
all tensors (\ref{y11}) invariant under $G$-action
(\ref{y10}).
Classification of all possible $W_{1234}$ leads to the following
explicit
form of bracket $(\ref{y12})$ (the detailed proof of this statement
will be published elsewhere):
\begin{equation}
\{\Omega_1,\Omega_2\}= [\Omega_{1} [ \Omega_{2}, \, r_{12}]]_+ +
X^{(1)}_{12}(\Omega_{1}^{2}+\Omega_{2}^{2})+
( \tilde{\Omega}_{1}^{2}+\tilde{\Omega}_{2}^{2})X^{(2)}_{12}+
\label{pb}
\end{equation}
$$
(\tilde{\Omega}_1X^{(3)}_{12}\Omega_1
+\tilde{\Omega}_2X^{(3)}_{12}\Omega_2 ) +
( \tilde{\Omega}_1X^{(4)}_{12}\Omega_2 +
\tilde{\Omega}_2X^{(4)}_{12}\Omega_1 ) +
X^{(5)}_{12}(\Omega_1\tilde{\Omega}_1+\Omega_2\tilde{\Omega}_2)+
$$
$$
( X^{(6)}_{12}(\tilde{\Omega}_1+\tilde{\Omega}_2)+(\Omega_1+\Omega_2)
X^{(7)}_{12} ) \mbox{tr}\Omega.
$$
where  all $X^{(i)}$ are symmetric $G$-invariant matrices in
$Mat(N,{\bf C})\otimes Mat(N,{\bf C})$:
$$
X^{(i)}= a_i I+b_i P+c_i K^{0}.
$$
and $a_i,~b_i,~c_i$ are complex numbers, $I$ is the identity matrix,
$P$-
is a permutation matrix and $K^{0}$:
$(K^{0})_{~ij}^{kl}=C^{kl}C_{ij}$.

Due to the identities
\begin{equation}
\label{y11a}
K^{0}_{12} \Omega_{1} = K^{0}_{12} \tilde{\Omega}_{2} \; , \;\;
K^{0}_{12} \Omega_{2} = K^{0}_{12} \tilde{\Omega}_{1} \; , \;\;
\end{equation}
we find that
$K^{0}_{12}(\Omega_1\tilde{\Omega}_1+\Omega_2\tilde{\Omega}_2)=0$,
i.e.,
we can put
$c_5=0$ and, therefore, the bracket (\ref{pb}) depends on
twenty arbitrary parameters $a_{i}, \; b_{i}, \; c_{i}$.
In fact this number coincides with dimension of
the cohomology group $H^{0}({\cal G},SV\otimes \wedge V)$, where
$V=Mat(N,{\bf C})$ and $SV~(\wedge V)$ stands for the symmetric
(antisymmetric) part of $V\otimes V$. We note that operators $X^{(i)}$
have the matrix structure of Yangian $R$-matrices.

Having the general form (\ref{pb}) one can calculate the bracket
between the variables $\Omega^{\pm}$. For this purpose one needs
explicit expressions for
$\{\tilde{\Omega}_1,\Omega_2\}$ and $\{{\Omega}_1,\tilde{\Omega}_2\}$
that are obtained from (\ref{pb}) by acting with $~\tilde{}~$ in the
first
or in the second matrix spaces. Now if we take into account
that
$$
\widetilde{X^{(i)}_{12}}(a,b,c)=X^{(i)}_{12}(a,\epsilon c,\epsilon b),
$$
then we get
\begin{equation}
\{\Omega_1^{\pm},\Omega_2^{\pm}\}=
[ \Omega_{1}^{\pm} [ \Omega_{2}^{\pm}, \, r_{12} ] ]_{+} +
Z_{12}^{\pm}(\Omega_{1}^{2}+\Omega_{2}^{2})-
(\tilde{\Omega}_{1}^{2}+\tilde{\Omega}_{2}^{2})Z^{\pm}_{12}+
\end{equation}
$$
( V_{12}^{\pm}(\Omega_1+\Omega_2)+(\tilde{\Omega}_{1}+
\tilde{\Omega}_{2})V_{12}^{\pm} ) \mbox{tr}\Omega,
$$
where
$$
Z_{12}^{\pm}=(X^{1}_{12}-X^{2}_{12})\pm (\tilde{X^{1}}-\tilde{X^{2}})=
\alpha^{\pm}Y_{12}^{\pm}+2\delta_{\pm+}(a_1-a_2),
$$
$$
V_{12}^{\pm}=(X^{6}_{12}+ X^{7}_{12})\pm
(\tilde{X^{6}}+\tilde{X^{7}})=
\beta^{\pm}Y_{12}^{\pm}+2\delta_{\pm+}(a_6+a_7),
$$
\begin{equation}
\begin{array}{l}
\alpha^{\pm}=b_1-b_2\pm\epsilon(c_1-c_2),\\
\beta^{\pm}=b_6+b_7\pm\epsilon(c_6+c_7)
\end{array}
\label{const}
\end{equation}
and $Y^{\pm}_{12}=P_{12}\pm \epsilon K^{0}_{12}$.
Thus, we see that Lie-valued
generators $\Omega^{-}$ form the closed algebra
\begin{equation}
\label{ppbb}
\{\Omega_1^-,\Omega_2^-\}= [\Omega_{1}^- [ \Omega_{2}^-, \,
r_{12}]]_+  \; ,
\end{equation}
only if $Z^{-}_{12} = V^{-}_{12} = 0$
or $\alpha^{-}=\beta^{-}=0$.
Then the calculation of the Jacobi identity (\ref{y5}) gives:
\begin{equation}
\label{jai}
\{\{ \Omega^{-}_1 , \, \Omega^{-}_{2} \}, \, \Omega^{-}_3 \} +
(cycle \; 1,2,3) =
- [ \Omega^{-}_{1}, \, [\Omega^{-}_{2}, \,
[\Omega^{-}_{3}, \, C(r)] ]_{+} ] \; ,
\end{equation}
where
\begin{equation}
\label{mybe}
C(r) = [r_{12}, \, r_{23} + r_{13}] + [r_{13}, \, r_{23} ] \; .
\end{equation}
If $r_{12}$ is a qusitriangular $r$-matrix
($C(r) \neq 0$ is $ad$-invariant tensor),
i.e., $r_{12}$ is a solution of the mYBE
(\ref{mybe}),
then the related bracket (\ref{ppbb})
is non Poisson (see (\ref{jai})).
Correspondingly, if $r_{12}$ is a triangular $r$-matrix ($C(r) = 0$
in eq. (\ref{mybe})),
then the bracket (\ref{ppbb}) is Poisson. These statements
agree with the results of \cite{AR}.

Before considering the general bracket (\ref{pb})
we recall how the exterior derivative $\bf d$ comes into this scheme.
If we relate in the quantum case
the co-invariant element $X$ of
Woronowicz with the quantum trace $\mbox{tr}_q\Omega$
(the definition of q-trace see in \cite{Fad},\cite{Resh},
\cite{Zum}, \cite{IsPop},
\cite{IsMal}), then
semiclassically it means that the ordinary exterior derivative $\bf
d$ is
expressed via the corresponding bicovariant bracket:
\begin{equation}
{\bf d}=\frac{1}{\kappa}\{ {\mbox
tr}\Omega,\ldots \},
\label{opd}
\end{equation}
where $\kappa$ is some
numerical parameter depending on a bracket under consideration.  The
fulfillment of the nilpotency condition: ${\bf d}^2=0$ is equivalent
to
the identity:
\begin{equation}
\{\{\Omega,\mbox{tr}\Omega\},\mbox{tr}\Omega\}=0;\\
\label{nil}
\end{equation}
and the Leibnitz rule is guaranteed by:
\begin{equation}
\{\{\Omega_1,\Omega_2\},\mbox{tr}\Omega\}+
\{\{\Omega_1,\mbox{tr}\Omega\},\Omega_2\}-
\{\Omega_1,\{\Omega_2,\mbox{tr}\Omega\}\}=0.
\label{dd}
\end{equation}
Bracket (\ref{pb}) satisfying (\ref{nil}) and (\ref{dd}) will be
refereed
as differential.
If (\ref{pb}) satisfy the Jacobi identity (\ref{y5}) then
(\ref{nil}) and (\ref{dd}) are fulfilled automatically.
It worth noting that $\mbox{tr}\Omega$
($(\mbox{tr}\Omega)^2=0$) looks like a BRST charge.

First, we find all differential brackets. From (\ref{pb})
one can extract the general forms
\begin{equation}
\{\Omega,\mbox{tr}\Omega\}=\mu_1\Omega^2+\mu_2\tilde{\Omega}^2+
\mu_3\tilde{\Omega}\Omega+\mu_4\Omega\tilde{\Omega}+
(\mu_5 \tilde{\Omega} + \mu_{6} \Omega)\mbox{tr}\Omega
\label{k2}
\end{equation}
and
\begin{equation}
\{\tilde{\Omega},\mbox{tr}\Omega\}=-\mu_1\tilde{\Omega}^2-
\mu_2\Omega^2
-\mu_3\tilde{\Omega}\Omega-\mu_4\Omega\tilde{\Omega}
+ (\mu_5 \Omega + \mu_{6} \tilde{\Omega})\mbox{tr}\Omega \; ,
\label{k3}
\end{equation}
where
\begin{equation}
\begin{array}{l}
\mu_1=2b_1+\epsilon c_1-\epsilon c_2+\epsilon c_4+N a_1,\\
\mu_2=-\epsilon c_1+\epsilon c_2+\epsilon c_4+2 b_2 +N a_2\\
\mu_3=\epsilon c_3+b _3+2 b_4 +N a_3,\\
\mu_4=\epsilon c_3-b_3+2 b_5 +N a_5,\\
\mu_5=a_4 + N a_{6} + 2 b_{6} + \epsilon (c_{6} + c_{7}), \\
\mu_6= - a_4 + N a_{7} + 2 b_{7} + \epsilon (c_{6} + c_{7}).
\end{array}
\label{k4}
\end{equation}
Substitution of (\ref{k4}) in (\ref{nil}) gives
four solutions for coefficients $\mu$:

i)
$\mu_1=\mu,~~~~ \mu_2=\mu_3=\mu_4=\nu,$ $~~\mu_{5} = \mu_{6} =0$, \\
where $\mu,\nu$ are arbitrary numbers except $\mu=\nu=0$;

ii) $\mu_1=\mu_2=-\mu_3=-\mu_4\neq 0$,
$\mu_{5} =\mu_{6} =0$;

iii) $\mu_5= -\mu_{6} = a_4 \neq 0$, $\mu_i=0 ~~(i=1,\dots, 4)$;

iv) $\mu_i=0 ~~(i=1,\dots, 6)$.

Thus, for bracket (\ref{k2}) (for (\ref{k3}) respectively) we have
four possibilities
\begin{equation}
\begin{array}{l}
i)~~~\{\Omega,\mbox{tr}\Omega\}=\mu\Omega^2+\nu(\tilde{\Omega}^2+
\tilde{\Omega}\Omega+\Omega\tilde{\Omega}),~~
(\mbox{for all}~ \mu,\nu \mbox{except}~\mu=\nu=0),\\
ii)~~\{\Omega,\mbox{tr}\Omega\}=\mu(\Omega^2+\tilde{\Omega}^2-
\tilde{\Omega}\Omega-\Omega\tilde{\Omega}),~~(\mu\neq 0),\\
iii)~\{\Omega,\mbox{tr}\Omega\}=\mu(\tilde{\Omega}-
\Omega)\mbox{tr}\Omega,~~
(\mu=a_4\neq 0),\\
iv)~~\{\Omega,\mbox{tr}\Omega\}=0.
\end{array}
\label{k5}
\end{equation}
The next step is to impose identity (\ref{dd}) proving the Leibnitz
rule
for the differential ${\bf d}$ (\ref{opd}).
It was done by using the symbolic manipulation program REDUCE. The
resulting differential bicovariant brackets are presented in Appendix.

Now substituting the calculated coefficients in (\ref{pb})
and analyzing the identity (\ref{y5}) with the help of the REDUCE
program, we arrive at the conclusion that neither of the
nontrivial differential brackets is Poisson.
Thus, among the family
(\ref{pb}) of bicovariant brackets there are
differential brackets but no Poisson brackets.
Note that we essentially use the requirement that
$\Omega$'s lie in the algebras 1.) $so(N)$, $sp(2n)$ or in
2.) $gl(N)=Mat(N)$. In the first case we have an additional
relation on the generators $\Omega^{+} = 0$. We stress that if we
consider
some other relations (cubic relations or
$\Omega_{1} \Omega_{2} K^{0}_{12} =
- K^{0}_{12}\Omega_{1} \Omega_{2}$),
then, the Poisson structure can exist.

Now we analyze
the external bicovariant algebra (\ref{w7}) on quantum groups
$SO_q(N)$ and
$Sp_q(N)$ directly in quantum case.
For these $q$-groups the $R$-matrix satisfies the cubic characteristic
equation \cite{Fad}:
\begin{equation}
\label{w18}
R = R^{-1} \, + \, \lambda \, - \, \lambda \,
K \; , \;\; K \equiv -\frac{1}{\lambda\nu} ( R^{2} -\lambda R -1) \; ,
\end{equation}
where
$\nu = \epsilon q^{\epsilon - N}$, $R=\hat{R}_{12} =P_{12}R_{12}$
and the matrix
$
K=K_{12} =
K^{i_{1}i_{2}}_{j_{1}j_{2}} = C^{i_{1}i_{2}}C_{j_{1}j_{2}}
$
is proportional to the singlet projector $P^{(0)}$:
\begin{equation}
\label{w19}
P^{(0)} = \mu^{-1} K \; , \;\;
\mu = (1 + \epsilon [N-\epsilon]_{q} ) \; .
\end{equation}
Note that this time $C$ is a quantum metric \cite{Fad}. Below we also
use the
projectors:
\begin{equation}
\label{w20}
P^{(\pm)} = \frac{1}{q+q^{-1}} (
\pm R + q^{\mp 1} I + \mu_{\pm} K) \; , \;\; \mu_{\pm} = -
\frac{q^{\mp 1}
\pm \nu}{\mu} \; .
\end{equation}
It has been shown in \cite{Wat}, \cite{Cast} that
for differential 1-forms one has the following relations coming from
definition (\ref{w7}):
\begin{equation}
\label{w19c}
X^{(\pm \pm)}=P^{(\pm)} \Omega' \, R \, \Omega' P^{(\pm)} = 0,~~~
X^{(00)}=P^{(0)} \Omega' \, R \, \Omega' P^{(0)} = 0 \; .
\end{equation}
Here $\Omega' = I \otimes \Omega = \Omega_{2}$
and the signs of the wedge products are omitted.
Taking the following sum
\begin{equation}
qX^{(++)}+\frac{1}{q}X^{(--)}-\frac{q\mu_{+}^2+q^{-1}\mu_{-
}^2}{(q+q^{-1})^2}
X^{(00)}=0
\label{weh}
\end{equation}
and using the identities
\begin{equation}
\label{ww}
\frac{\mu_{+} + \mu_{-}}{q + q^{-1}} = - \frac{1}{\mu} \; , \;\;
\frac{q\mu_{+} - q^{-1}\mu_{-}}{q + q^{-1}} = -\frac{\nu}{\mu} \; ,
\end{equation}
one can show that relations (\ref{w19c}) are
equivalent to the unique relation:
\begin{equation}
\label{w22}
(R \, \Omega' \, R \, \Omega' \, R + \Omega' \, R  \, \Omega')
- \frac{1}{\mu}
( K \, \Omega' R \, \Omega' + \Omega' \, R \, \Omega' \, K) -
\end{equation}
$$
- \frac{\nu}{\mu}
( K \, \Omega' \, R \, \Omega' \, R +
R \, \Omega' \, R \, \Omega' \, K)  = 0 \; .
$$
This form of the defining relations for $\Gamma^{\wedge}$ is suitable
to produce a graded bicovariant bracket on $\Gamma^{\wedge}_{cl}$.

The semiclassical expansions of projector $P^{(0)}_{12}$
and $R$-matrix are:
$$
P^{(0)}_{12}=\hat{P}^{(0)}_{12}+
\hbar\frac{\epsilon}{N}K_{12}^{1}+O(\hbar^2),
{}~~~R_{12}=P_{12}+\hbar P_{12}\tilde{r}_{12}+O(h^2),
$$
where $\hat{P}^{(0)}_{12}=\frac{\epsilon}{N}K_{12}^{0}$ and
$\tilde{r}$ satisfies CYBE.
It follows from $KR=RK=\nu K$ that in the first order in $\hbar$:
\begin{equation}
\begin{array}{l}
K_{12}^{1}-\epsilon K_{12}^{1}P_{12}=K_{12}^{0}\tilde{r}_{12}-
\epsilon(1-\epsilon N)
K_{12}^{0},\\
K_{12}^{1}-\epsilon P_{12} K_{12}^{1}=\tilde{r}_{21}
K_{12}^{0}-\epsilon(1-\epsilon N) K_{12}^{0}.
\end{array}
\label{qe}
\end{equation}
Then, by expanding (\ref{w22}) in
powers of $\hbar$,
taking into account (\ref{qe}) and the correspondence
principle (\ref{cl}),  we obtain:
\begin{equation}
(I-\hat{P}^{(0)}_{12})(\{\Omega_1,\Omega_2\}+G_{12})-
(\{\Omega_1,\Omega_2\}+G_{12})\hat{P}^{(0)}_{12}=0,
\label{genw}
\end{equation}
where
\begin{equation}
G_{12}=-[\Omega_1,[\Omega_2,r_{12}]]_{+}
+P_{12}(\Omega_1^2+\Omega_2^2)-
\epsilon(K_{12}\Omega_1\Omega_2+\Omega_1\Omega_2K_{12}+
\Omega_1K_{12}\Omega_2+\Omega_2K_{12}\Omega_1)
\label{gg}
\end{equation}
and we made use of
the quasitriangular $r$-matrix: $r=\tilde{r}-(P-\epsilon K)$.

The components $\hat{P}^{(0)}_{12} \{ \Omega_{1}, \, \Omega_{2} \}
(I-\hat{P}^{(0)}_{12})$ are not defind by (\ref{genw}).
Thus, we see that relations (\ref{w19c}) are unsufficient to
generate, in the limit $\hbar\rightarrow 0$, a genuine bicovariant
bracket.
In the quantum case it means that
the number of defining relations (\ref{w19c}) is not enough
to reoder lexicographically arbitrary monom in $\Omega$'s.
Therefore, if we confine ourselves only with (\ref{w19c}),
then we can not conclude that $\dim \Gamma^{\wedge}$ is equal to
$\dim \Gamma^{\wedge}_{cl}$.

On the other hand, we can not assume the solution of
(\ref{genw}) as
\begin{equation}
\label{Poi}
\{ \Omega_{1}, \, \Omega_{2} \} = - G_{12} \; ,
\end{equation}
since $G_{12}$ is symmetric under $1 \leftrightarrow 2$
only if the following relation hold:
\begin{equation}
\label{gru}
K_{12} \Omega_{1} \Omega_{2} +
\Omega_{1} \Omega_{2} K_{12} = 0 \; .
\end{equation}
But this relation contradicts to the requirement that
$\Omega \in {\cal G}$ or that the number of $\Omega$s are $N^{2}$
($\Omega \in Mat(N)$).
Note, however, that the bracket (\ref{Poi}) is Poisson
for $\Omega$s restricted by constraint (\ref{gru}).

To improve the situation
the authors of \cite{Wat}, in addition to (\ref{w19c}),
have assumed the relations:
\begin{equation}
X^{(0+)}=P^{(0)} \Omega' R \Omega'
P^{(+)}=0,~~~X^{(+0)}=P^{(+)} \Omega' R \Omega' P^{(0)} = 0.
\label{rel1}
\end{equation}
One can obtain without problems that (\ref{w22}) and (\ref{rel1})
are equivalent to the relation:
\begin{equation}
R \, \Omega' \, R \, \Omega' \, R + \Omega' \,
R  \, \Omega' +\frac{1}{\mu}(\nu q^{-1}-1) ( K \, \Omega' R \,
\Omega' +
\Omega' \, R \, \Omega' \, K)=0.
\label{wat}
\end{equation}
By expanding
(\ref{wat}) in $\hbar$, as it was done for the general
relation (\ref{w22}), we get the
following bicovariant bracket:
\begin{equation}
\{\Omega_1,\Omega_2\}=[\Omega_1,[\Omega_2,r_{12}]]_+
-P_{12}(\Omega_1^2+\Omega_2^2)+
(\Omega_1K_{12}\Omega_2+\Omega_2K_{12}\Omega_1).
\label{fgf}
\end{equation}
being a particular case of (\ref{pb}). Now we see
that according to our classification this bracket neither
Poisson nor differential.
It means in the quantum case that the requirement ${\bf d}^2=0$
(\ref{nil}) implies some additional cubic
relations on generators $\Omega^{i}_{j}$,
which were not assumed from the beginning.
The situation is somewhat improved when we require ${\bf d}^2=0$
only on "physical" components $\Omega =\Omega^-$.
This requirement is consistent with (\ref{nil}),
since $\{\Omega^-,\mbox{tr}\Omega\}=-2(\Omega^-)^{2}$.
However, the substitution
$\Omega \rightarrow \Omega^-$ in (\ref{dd}) leads to the
conclusion:
$$
\{\{\Omega^{-}_1,\Omega^{-}_2\},\mbox{tr}\Omega\}+
\{\{\Omega^{-}_1,\mbox{tr}\Omega\},\Omega^{-}_2\}-
\{\Omega^{-}_1,\{\Omega^{-}_2,\mbox{tr}\Omega\}\} \neq 0.
$$
Thus, one can not assume the
Leibnitz rule for $\bf d$ on the
"physical" subalgebra generated by $\Omega^-$ without
imposing new cubic relations on $\Omega$'s. Note that
if we impose the unacceptable relations (\ref{gru}), then,
the bracket (\ref{fgf}) coinsides with (\ref{Poi}) and, therefore,
is Poisson.

Seemingly, the absence of bicovariant
Poisson structure for $SO(N)$ and $Sp(N)$ ($N$ is generic) reflects
the fact that we can not confine ourselves by considering
only $G$-invariant tensors $W$ in (\ref{y12}).
Considering in (\ref{y12})
tensor $W$ which is not $G$-invariant, we disturb, of course,
bicovariance but may hope to keep the Jacobi identity.
Then we expect that the bicovariance will be restored
on the surface $\Omega^+=0$
if we treat $\Omega^+=0$ as the
first order constraint (in the Dirac sense).

$$~$$ {\bf ACKNOWLEDGMENT} $$~$$
The authors are grateful to J.Lukierski for interesting discussions
and to P.B.Medvedev and P.N.Pyatov for helpful comments.
This work was partially supported by the grant KBN 2P 30208706,
RFFR (grants N93-011-147, N93-02-3827),  by ISF (grants
M1L-000, RFF 000) and INTAS (grant 93-127).

{\large \bf APPENDIX}
\appendix
\section{Differential $SO$- and $Sp$- covariant
brackets on $Mat(N)$.}

1. First solution:
$\{\Omega,\mbox{tr}\Omega\}=\mu\Omega^2+\nu(\tilde{\Omega}^2+
\tilde{\Omega}\Omega+\Omega\tilde{\Omega})$

i)
$$
\{\Omega_1,\Omega_2\}=[\Omega_1,[\Omega_2,r_{12}]]_+  -
\frac{1}{N}(2b_2+\epsilon c_3+Na_2)(\Omega_1^2+\Omega_2^2)
$$
$$
+a_2((\Omega^+_1)^2+(\Omega^+_2)^2)+c_3\Omega^{+}_{1}K_{12}\Omega^{+}_
\frac{1}{2}\mu P_{12}(\Omega_1^2+\Omega_2^2)
$$
$$
+b_2P_{12}((\Omega^+_1)^2+(\Omega^+_2)^2)+
c_6(K_{12}\Omega^+_1-\Omega^+_1K_{12})\mbox{tr}\Omega;
$$

ii)
$$
\{\Omega_1,\Omega_2\}=[\Omega_1,[\Omega_2,r_{12}]]_+  -
\frac{1}{N}(2b_1+\epsilon c_3+Na_2)(\Omega_1^2+\Omega_2^2)
$$
$$
+a_2((\Omega^+_1)^2+(\Omega^+_2)^2)+c_3\Omega^{+}_{1}K_{12}\Omega^{+}_
(2b_1+\nu)P_{12}((\Omega^+_1)^2+(\Omega^+_2)^2)+
$$
$$
-(b_1+\nu)(P_{12}(\Omega_1^2+\Omega_2^2)+
(\tilde{\Omega}_1^2+\tilde{\Omega}_2^2)P_{12})
-\frac{1}{2}\epsilon \nu
(K_{12}(\Omega_1^2+\Omega_2^2)+(\Omega_1^2+\Omega_2^2)K_{12})+
$$
$$
c_6(K_{12}\Omega^+_1-\Omega^+_1K_{12})\mbox{tr}\Omega;
$$

iii)
$$
\{\Omega_1,\Omega_2\}=[\Omega_1,[\Omega_2,r_{12}]]_+  -
\frac{1}{N}(2b_2+\epsilon c_3-\nu)((\Omega^+_1)^2+(\Omega^+_2)^2)
$$
$$
+c_3\Omega^{+}_{1}K_{12}\Omega^{+}_{2}+
b_1 P_{12}((\Omega^+_1)^2+(\Omega^+_2)^2)+
$$
$$
(a_6(\tilde{\Omega}_1+\tilde{\Omega_2})-
3a_6(\Omega_1+\Omega_2)
-\frac{1}{2}\epsilon(c_6+c_7)P_{12}
(\Omega^+_1+\Omega^+_2)+c_6K_{12}\Omega^+_1+c_7\Omega^+_1K_{12})
\mbox{tr}\Omega;
$$

2. Second solution: $\{\Omega,\mbox{tr}\Omega\}=
\mu(\Omega^2+\tilde{\Omega}^2-
\tilde{\Omega}\Omega-\Omega\tilde{\Omega})$.

$$
\{\Omega_1,\Omega_2\}=[\Omega_1,[\Omega_2,r_{12}]]_+  +
$$
$$
a_1((\Omega^-_1)^2+(\Omega^-_2)^2)+c_3\Omega^{-}_{1}K_{12}\Omega^{-
}_{2}
-\frac{1}{2}b_3P_{12}(\Omega_1^+\Omega_1^- + \Omega_2^+\Omega_2^-)+
$$
$$
(-\epsilon c_1 P_{12}+c_1 K_{12})(\Omega_1^2+\Omega_2^2)+
(\tilde{\Omega}_1^2+\tilde{\Omega}_2^2)
(\epsilon c_1 P_{12}+c_2K_{12})
$$
$$
((a_6+\epsilon c_6P_{12}+c_6 K_{12})\Omega_1^++
\Omega_1^+(a_6+\epsilon c_6P_{12}+c_6 K_{12}))\mbox{tr}\Omega,
$$
where $b_3=-\epsilon(c_1+c_2)$.

3.Third solution: $\{\Omega,\mbox{tr}\Omega\}=\mu(\tilde{\Omega}-
\Omega)\mbox{tr}
\Omega$.
$$
\{\Omega_1,\Omega_2\}=[\Omega_1,[\Omega_2,r_{12}]]_+
$$
$$
+(-a_3-b_4+c_1K_{12})(\Omega_1^2+\Omega_2^2)+
(\tilde{\Omega}_1^2+\tilde{\Omega}_2^2)(a_3+b_4+c_1K_{12})
$$
$$
(a_3+b_4 P_{12})(\Omega_1^+\Omega_1^- + \Omega_2^+\Omega_2^-)+
(\mu+b_3P_{12})(\tilde{\Omega}_2\Omega_1+\tilde{\Omega}_1\Omega_2)+
$$
$$
(X_{12}^{(6)}(\tilde{\Omega}_1+\tilde{\Omega_2})+
(\Omega_1+\Omega_2)X_{12}^{(7)})\mbox{tr}\Omega,
$$
where $b_3=-(Na_3+2b_4)$ and coefficients in $X_{12}^{(6)}$ and
in $X_{12}^{(7)}$ remain to be arbitrary.
\end{document}